\documentclass[a4paper,10pt,oneside]{amsart} 

\usepackage{a4wide}
\usepackage[latin1]{inputenc}
\usepackage{amssymb}
\usepackage{url}
\usepackage[bookmarks=false,pdfborder={0 0 0.05}]{hyperref}
\usepackage{color}
\usepackage{tikz}
\usetikzlibrary{calc,through,backgrounds}

\newtheorem{theorem}{Theorem}[section]
\newtheorem{proposition}[theorem]{Proposition}
\newtheorem{corollary}[theorem]{Corollary}
\newtheorem{lemma}[theorem]{Lemma}

\theoremstyle{definition}

\newtheorem{example}[theorem]{Example}
\newtheorem{examplecon}[theorem]{Example (continued)}
\newtheorem*{remark}{Remark}  

\def\Sn{\mathcal{S}_n}

\def\defn#1{{\sf #1}}

\def\D{\mathcal{D}}

\def\maj{\operatorname{maj}}
\def\fmaj{\operatorname{fmaj}}
\def\imaj{\operatorname{imaj}}
\def\ifmaj{\operatorname{ifmaj}}
\def\Des{\operatorname{Des}}
\def\iDes{\operatorname{iDes}}
\def\fDes{\operatorname{fDes}}
\def\ifDes{\operatorname{ifDes}}
\def\des{\operatorname{des}}
\def\ides{\operatorname{ides}}
\def\fdes{\operatorname{fdes}}
\def\ifdes{\operatorname{ifdes}}
\def\Asc{\operatorname{Asc}}

\def\Cat{\operatorname{Cat}}
\def\aCat{\operatorname{Cat}}

\def\area{\operatorname{area}}
\def\Set{\operatorname{Set}}
\def\inv{\operatorname{inv}}
\def\rev{\operatorname{rev}}
\def\neg{\operatorname{neg}}
\def\Neg{\operatorname{Neg}}

\def\NC{N\!C}
\def\Dn{\mathcal{D}_n}

\def\Cox{\operatorname{Sort}}
\def\ZZ{\mathbb{Z}}

\definecolor{grey}{rgb}{.5 , .5 , .5}
\definecolor{lightgrey}{rgb}{.9 , .9 , .9}

\newcommand{\latticepath}[4]{ 
  \coordinate (L) at #2;
  \foreach \x/\y/\a in {#4} {
    \coordinate (L1) at ($ (L) + ( #1 * \x , #1 * \y ) $);
    \draw[color=\a, #3] (L) -- (L1);
    \coordinate (L) at (L1);
  }
}

\newcommand{\tikzbox}[5]{ 
  \coordinate (A) at #2;
  \coordinate (B) at ($ (A) + (  #1,-#1  ) $);
  \coordinate (C) at ($ (A) + (#1/2,-#1/2) $);
  \draw[#3, very thin, fill=#4] (A) rectangle (B);
  \node at (C) []{#5};
}

\newcommand{\greybox}[3]{ 
  \tikzbox{#1}{#2}{grey}{white}{#3}
}

\newcommand{\boxcollection}[3]{ 
  \coordinate (X) at #2;
  \foreach \x/\y/\object in {#3} {
    \greybox{#1}{($ (X) + (#1 * \y, - #1 * \x) - (#1, - #1) $)}{\object};
  }
}

\newcommand{\content}[3]{ 
  \coordinate (C) at #2;
  \foreach \x/\y/\object in {#3} {
    \node at ($ (C) + ( #1 * \y, -#1 * \x ) - ( #1 / 2, -#1 / 2 )$) {\object};
  }
}

\begin{document}
  \title[More bijective Catalan combinatorics on permutations]{More bijective Catalan combinatorics on\\permutations and on signed permutations}
  \author{Christian Stump}
  \address{Freie Universit\"at Berlin, Germany}
  \email{christian.stump@fu-berlin.de}
  \urladdr{http://homepage.univie.ac.at/christian.stump/}
  \date{\today}

  \begin{abstract}
    In this paper, we construct bijections between Dyck paths, noncrossing partitions, and $231$-avoiding permutations, which send the area statistic on Dyck paths to the inversion number on noncrossing partitions and on $231$-avoiding permutations. This bijection has the additional property that it simultaneously sends the major index on Dyck paths to the sum of the major index and the inverse major index on noncrossing partitions and on $231$-avoiding permutations, respectively. Moreover, we provide generalizations of these constructions to the group of signed permutations.
  \end{abstract}

  \maketitle

  \tableofcontents
  
\section{Introduction}

  A set partition $P$ of a totally ordered finite set $S$ is called \defn{noncrossing} if there do not exist elements $a<b<c<d$ in $S$ such that $a$ and $c$ are in one block of $P$ while $b$ and $d$ are in another block.
  Graphically, this means that when the elements of $S$ are drawn in the given order on a circle, that $P$ is noncrossing if and only if the convex hulls of the blocks of $P$ do not cross.
  In particular, the property of being noncrossing only depends on the cyclic ordering of~$S$.
  One can naturally embed all set partitions of $\{1,\ldots,n\}$ into the set $\Sn$ of permutations of $\{1,\ldots,n\}$ by sending a set partition $P$ to the permutation whose cycles are the blocks of $P$ in increasing order.
  E.g., the noncrossing set partition $\big\{ \{ 1,3,9\},\{2\},\{4\},\{5,6,7,8\} \big\}$ of $\{1,2,\ldots,9\}$ is sent to the permutation $(1,3,9)(5,6,7,8) \in \mathcal{S}_9$, written here in cycle notation.
  We denote the set of noncrossing partitions, considered inside $\Sn$, by $\NC_n$.

  A permutation $\sigma \in \Sn$ is called \defn{$231$-avoiding} if there do not exist $a<b<c$ such that $\sigma_c < \sigma_a < \sigma_b$, where we write $\sigma = [\sigma_1,\ldots,\sigma_n]$ in one-line notation.
  In other words,~$\sigma$ is $231$-avoiding if $[\sigma_1,\ldots,\sigma_n]$ does not contain a subword that is order-isomorphic to $231$.
  We denote the set of $231$-avoiding permutations in $\Sn$ by $\Cox_n$. This notation refers to the fact that being $231$-avoiding is equivalent to being \emph{stack sortable}, see~\cite{Rea2007}.

  It is well-known that noncrossing partitions and $231$-avoiding permutations are both counted by the Catalan numbers,
  $$|\NC_n| = |\Cox_n| = \Cat_n := \frac{1}{n+1}\binom{2n}{n}.$$

  \smallskip

  In this paper we describe connections between several statistics on noncrossing partitions, $231$-avoiding permutations, and another interesting family counted by the Catalan numbers, namely \defn{Dyck paths}.
  These are lattice paths in $\ZZ^2$ starting at $(0,0)$, consisting of~$n$ north and~$n$ east steps and which never go below the diagonal $x=y$.
  We denote the set of all Dyck paths of length~$2n$ by $\Dn$.

\subsection{Results for the group of permutations}

  The first theorem concerns the \defn{major index} of $\sigma \in \Sn$ defined as
  $$\maj(\sigma) = \sum_{i \in \Des(\sigma)}i,$$
  and the \defn{inverse major index} defined as $\imaj(\sigma) = \maj(\sigma^{-1})$. Here, the \defn{descent set} of a permutation~$\sigma$ is given by $\Des(\sigma) = \{ 1 \leq i < n: \sigma_i>\sigma_{i+1} \}$. For later convenience, we define also $\iDes(\sigma) = \Des(\sigma^{-1})$, $\des(\sigma) = |\Des(\sigma)|$, and $\ides(\sigma) = |\iDes(\sigma)|$.
  Moreover, define the inversion number as
  $$\inv(\sigma) = \big| \{ 1 \leq i < j \leq n : \sigma_i > \sigma_j \} \big|.$$
  Finally, let $[k]_q = 1+q+\ldots+q^{k-1}$ denote the usual $q$-extension of the integer $k$, $[k]_q! = [1]_q[2]_q \cdots [k]_q$ be the $q$-factorial, and $\left[\begin{smallmatrix} k \\ l \end{smallmatrix}\right]_q = [k]_q!/[l]_q![k-l]_q!$ the $q$-binomial coefficient.

  \begin{theorem}\label{th:majimajNCCox}
    The generating function for $\maj + \imaj$ on noncrossing partitions and on $231$-avoiding permutations is given by MacMahon's $q$-Catalan numbers,
    \begin{align*}
      \sum_{\sigma \in \NC_n} q^{\maj(\sigma) + \imaj(\sigma)} = \sum_{\sigma \in \Cox_n} q^{\maj(\sigma) + \imaj(\sigma)} = \frac{1}{[n+1]_q}\begin{bmatrix} 2n \\ n \end{bmatrix}_q.
    \end{align*}
  \end{theorem}

  We prove this theorem by providing bijections between noncrossing partitions and $231$-avoiding permutations in $\Sn$ on the one hand and Dyck paths of length~$2n$ on the other hand.
  For a Dyck path $D$, the \defn{major index} is given by
  $$\maj(D) := \sum_{i \in \Des(D)}i,$$
  where $\Des(D)$ is given by the collection of indices~$i$ such that the~$i$-th step in $D$ is an east step and the $(i+1)$-st step is a north step.
  In \cite{Mah1960}, P.A.~MacMahon showed that its generating function is given by the above expression,
  $$\sum_{D \in \Dn} q^{\maj(D)} = \frac{1}{[n+1]_q}\begin{bmatrix} 2n \\ n \end{bmatrix}_q.$$
  Another natural statistic on Dyck paths is the \defn{area statistic}. For $D \in \Dn$, the statistic $\area(D)$ is defined to be the number of open boxes
  \begin{align}
    b_{ij} = \{ (x,y) \in \mathbb{R}^2 : i < x < i+1, j < y < j + 1 \} \label{eq:boxes}
  \end{align}
  which lie below $D$ and strictly above the diagonal $x=y$. The following theorem implies Theorem~\ref{th:majimajNCCox}, as the coefficients of MacMahon's $q$-Catalan numbers are known to be symmetric of degree $n(n-1)$.
  \begin{theorem}\label{th:bijectionsA}
    There exist explicit bijections $\phi_n : \Dn \tilde\longrightarrow \NC_n$ and $\psi_n : \Dn \tilde\longrightarrow \Cox_n$ such that for all $D \in \Dn$,
    \begin{align}
      \area(D)            &= \inv(\phi_n(D)) = \inv(\psi_n(D)), \label{eq:areatoinv} \\
      n(n-1) - \maj(D)    &= \maj(\phi_n(D)) + \imaj(\phi_n(D)) = \maj(\psi_n(D)) + \imaj(\psi_n(D)) \label{eq:majtomaj}.
    \end{align}
  \end{theorem}

  \subsection{Extensions to the group of signed permutations}

  The various objects and statistics in the previous section have natural analogues for the group $B_n$ of signed permutations. A \defn{signed permutation} $\sigma \in B_n$ is a bijection of the integers $\{ \pm 1, \ldots, \pm n \}$ such that $\sigma_{-i} = -\sigma_i$. Thus,~$\sigma$ can be represented in \defn{one-line notation} by $\sigma = [ \sigma_1,\ldots,\sigma_n]$. We will also use the \defn{cycle notation} for signed permutations. For example, we write
  $$[-3,-2,1] = (1,-3,-1,3)(2,-2) \in B_3.$$
  The notation $B_n$ is due to the fact that the group of signed permutations is the Weyl group of Cartan-Killing type $B_n$.
  We now state the extended theorems for signed permutations. For definitions, in particular for the definition of the flag major index and of type $B_n$ reverse noncrossing partitions, we refer to Section~\ref{sec:signedpermutations}.
  \begin{theorem}\label{th:bijectionsB}
    There exist explicit bijections $\phi_{B_n} : \D_{B_n} \tilde\longrightarrow \rev(\NC_{B_n})$ and $\psi_{B_n} : \D_{B_n} \tilde\longrightarrow \Cox_{B_n}$ such that for $D \in \D_{B_n}$,
    \begin{align}
      \area(D) &= \inv_B(\phi_{B_n}(D)) = \inv_B(\psi_{B_n}(D)), \label{eq:areatoinvB} \\
      2n^2 - \fmaj(D) &= \fmaj(\phi_{B_n}(D)) + \ifmaj(\phi_{B_n}(D)) = \fmaj(\psi_{B_n}(D)) + \ifmaj(\psi_{B_n}(D)). \label{eq:majtomajB}
    \end{align}
  \end{theorem}
  As for permutations, this implies the following corollary. Here, we take again into account that the coefficients of MacMahon's $q$-Catalan numbers for $B_n$, as given in the corollary below, are symmetric, and that its degree is $2n^2$.
  \begin{corollary}\label{cor:majimajB}
    The generating function for $\fmaj+\ifmaj$ on reverse noncrossing partitions and on sortable elements for signed permutations are given by MacMahon's $q$-Catalan numbers for $B_n$,
    \begin{align*}
      \sum_{\sigma \in \rev(\NC_{B_n})} q^{\fmaj(\sigma) + \ifmaj(\sigma)} = \sum_{\sigma \in \Cox_{B_n}}  q^{\fmaj(\sigma) + \ifmaj(\sigma)} =  \begin{bmatrix} 2n \\ n \end{bmatrix}_{q^2}.
    \end{align*}
  \end{corollary}

  \begin{remark}
    All objects and statistics in the previous corollary have as well analogues for the group of \emph{even-signed permutations}, which is the Weyl group of Cartan-Killing type $D_n$. Analogous statements do not hold in that case.
  \end{remark}

\section{Proofs for the group of permutations}\label{sec:permutations}

    In this section we construct the proposed bijections $\phi_n : \Dn \tilde\longrightarrow \NC_n$ and $\psi_n : \Dn \tilde\longrightarrow \Cox_n$ and show that they have the desired properties.

    \subsection{Dyck paths and noncrossing partitions}
\begin{figure}
          \centering
          \begin{tikzpicture}[scale=1]
            \boxcollection{0.5}{(0,0)}{
            1/6/,1/7/,1/8/,
            2/5/,2/6/,2/7/,
            3/4/,3/5/,3/6/,
            4/4/,4/5/,
            5/3/,5/4/,
            6/3/,
            7/1/,7/2/,
            8/1/}
            \latticepath{0.5}{(0,-4.5)}{very thick}{
              0/3/black,2/0/black,0/2/black,1/0/black,0/2/black,1/0/black,0/1/black,1/0/black,0/1/black,4/0/black}
            \latticepath{0.5}{(0.25,-4.25)}{very thick}{
              0/2/grey,2/0/grey,0/2/grey,1/0/grey,0/2/grey,1/0/grey,0/1/grey,1/0/grey,0/1/grey,3/0/grey}
            \latticepath{0.5}{(2.25,-2.25)}{very thick}{
              0/1/grey,1/0/grey,0/1/grey,1/0/grey,0/1/grey,1/0/grey}
            \latticepath{0.5}{(0,-4.5)}{thin}{
              9/9/black}
            \content{0.5}{(0.5,-4.5)}{
              0/0/\colorbox{white}{$\scriptstyle 1$},-1/1/\colorbox{white}{$\scriptstyle 2$},-2/2/\colorbox{white}{$\scriptstyle 3$},-3/3/\colorbox{white}{$\scriptstyle 4$},-4/4/\colorbox{white}{$\scriptstyle 5$},-5/5/\colorbox{white}{$\scriptstyle 6$},-6/6/\colorbox{white}{$\scriptstyle 7$},-7/7/\colorbox{white}{$\scriptstyle 8$},-8/8/\colorbox{white}{$\scriptstyle 9$}}
          \end{tikzpicture}
          \caption{The bijection $\phi_9$ sending the shown Dyck path to the noncrossing partition $\sigma = (1,3,9)(5,6,7,8) = [3,2,9,4,6,7,8,5,1]$.}
          \label{fig:bijNNNCA}
        \end{figure}
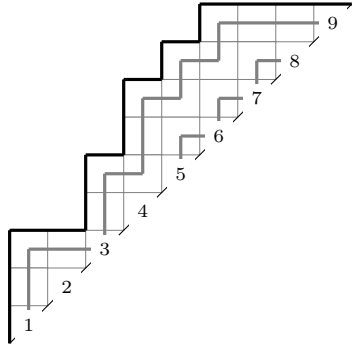

      Write the numbers $1$ through~$n$ on the diagonal below the Dyck path $D$ by labeling the box $b_{ii}$ by $i+1$ for $0 \leq i < n$.
      Define $\phi_n(D)$ to be the permutation obtained by a \lq\lq shelling\rq\rq\ of $D$ and connecting the appropriate numbers as indicated in Figure~\ref{fig:bijNNNCA}. For the Dyck path $NNNEENNENNENENEEEE \in \D_9$, the \lq\lq shell\rq\rq\ connecting $1,3$ and $9$ forces $(1,3,9)$ to form an increasing cycle in $\phi_n(D)$. After shelling the second layer with the path connecting $5,6,7$, and $8$, the complete Dyck paths is \lq\lq shelled\rq\rq. Thus, we have $\phi_9(D) = (1,3,9)(5,6,7,8)$.

      The following proposition is Theorem~\ref{th:bijectionsA}\eqref{eq:areatoinv} for $\phi_n : \Dn \tilde\longrightarrow \NC_n$.
      \begin{proposition} \label{prop:area}
       ~$\phi_n$ is a bijection between $\Dn$ and $\NC_n$ and maps the area statistic on $\Dn$ to the inversion number on $\NC_n$. For $D\in\Dn$,
        $$\area(D) = \inv(\phi_n(D)).$$
      \end{proposition}
      \begin{proof}
        The described construction is clearly well-defined and invertible, and thus,~$\phi_n$ is a bijection. To prove that the area statistic is mapped to the inversion number, we assume first that $\phi_n(D)$ consists of a unique nontrivial cycle, $\phi_n(D) = (i_1,\ldots,i_k)$ with $k>1$. The general case is obtained by applying the same argument several times. Observe that we can indeed analyze each cycle individually as $\phi_n(D)$ is noncrossing and therefore, the inversion number of one cycle of $\phi_n(D)$ is independent of the other cycles, and that the area in one shell of $D$ only depends on its length and the number of times it passes through the cells on the main diagonal.
        The area of $D$ is equal to $2(i_k-i_1)-1-(k-2) = 2(i_k-i_1)-k+1$. It is easy to see that this is also equal to the inversion number of the cycle $(i_1,\ldots,i_k)$.
      \end{proof}
      \begin{example}
        Let $D$ and $\sigma = \phi_9(D)$ as in Figure~\ref{fig:bijNNNCA}. Then the \lq\lq shell\rq\rq\ connecting $1,3,$ and $9$ contains $14$ boxes which is equal to the inversion number of the cycle $(1,3,9)$, the \lq\lq shell\rq\rq\ connecting $5,6,7$, and $8$ contains $3$ boxes which is equal to the inversion number of the cycle $(5,6,7,8)$.
      \end{example}
      The following proposition is Theorem~\ref{th:bijectionsA}\eqref{eq:majtomaj} for $\phi_n : \Dn \tilde\longrightarrow\NC_n$.
      \begin{proposition}\label{prop:majimajA}
        Let $D \in \Dn$. Then
        $$\maj(D) + \maj(\phi_n(D)) + \imaj(\phi_n(D)) = n(n-1).$$
      \end{proposition}
      Before proving the theorem in several steps, we get back to the example.
      \begin{examplecon}
        The descent set of the Dyck path $D$ in Figure~\ref{fig:bijNNNCA} and the descent set and the inverse descent set of $\sigma = \phi_9(D)$ are given by
        $$\Des(D) = \{5,8,11,13\}, \quad \Des(\sigma) = \{1,3,7,8\}, \quad \iDes(\sigma) = \{1,2,5,8\}$$
        and therefore,
        \begin{align*}
          \maj(D) &= 5+8+11+13 = 37, \\
          \maj(\sigma)+\imaj(\sigma) &= (1+3+7+8) + (1+2+5+8) = 19 + 16, \\
          \maj(D) + \maj(\sigma) + \imaj(\sigma) &= 37 + 19 + 16 = 72 = 9 \cdot 8.
        \end{align*}
      \end{examplecon}      
      \begin{lemma}\label{le:NNdes}
        Let $\sigma \in \NC_n$.
        Then $\des(\sigma) = \ides(\sigma)$.
      \end{lemma}
      \begin{proof}
        We first prove the lemma for the case that~$\sigma$ has a unique nontrivial cycle $\sigma = (i_1,\ldots,i_k)$. As~$\sigma \in \NC_n$, we can assume that $i_1 < \ldots < i_k$. We can thus explicitly describe the descent set and the inverse descent set of~$\sigma$,
        \begin{align*}
          \Des(\sigma)  &= \{ i_\ell : \ell < k, i_\ell+1 < i_{\ell+1}\} \cup \{i_k-1\},\\
          \iDes(\sigma) &= \{i_1\} \cup \{ i_\ell-1 : 1 < \ell, i_{\ell-1} + 1 < i_\ell \}.
        \end{align*}
        Obviously, both have the same size.
        Next, observe that the descent and inverse descent sets of any $\sigma \in \NC_n$ is given by the disjoint union of the descent and inverse descent sets of each of the cycles of~$\sigma$ individually.
        This comes from the fact that the cycles of~$\sigma$ do not cross which yields the observation that the descent and inverse descent sets of~$\sigma$ do not interfere.
        This completes the proof of the lemma.
      \end{proof}
      We will use the description of the descent set and the inverse descent set given in the previous proof as well in subsequent constructions. Since this description depends on blocks of consecutive integers in cycles of~$\sigma$, we define a \defn{block} of~$\sigma$ to be  a maximal part $i_\ell,i_{\ell+1}=i_\ell+1,\ldots,i_{\ell+a}=i_\ell+a$ of consecutive integers within a nontrivial cycle $(i_1,\ldots,i_k)$ of~$\sigma$. E.g., $(1\ |\ 3\ |\ 9)(5,6,7,8)$ in Figure~\ref{fig:bijNNNCA} has $4$ blocks, namely $1, 3, 9,$ and $5678$. As in this example, we sometimes separate blocks by vertical bars.
      \begin{lemma} \label{le:descentvalleys}
        Let $D \in \Dn$ and let $\sigma = \phi_n(D)$. Then
        $$\des(D) + \des(\sigma) = n-1.$$
      \end{lemma}
      \begin{proof}
        From the description of the descent set of~$\sigma$ in the proof of Lemma~\ref{le:NNdes}, we see that $\des(\sigma)$ equals the number of blocks of~$\sigma$. On the other hand, every peak of $D$ corresponds either to a fixpoint of $\phi_n(D)$ if it lies on an anti-diagonal above an integer on the main diagonal --~this happens for $2$ and $4$ in Figure~\ref{fig:bijNNNCA}, which are exactly the fixpoints of $\sigma = (1,3,9)(5,6,7,8)$~-- or it corresponds to two consecutive integers in a block if it lies on an anti-diagonal between these two integers on the main diagonal --~this happens for $56,67$, and $78$ in Figure~\ref{fig:bijNNNCA}, which are exactly the connections in the big block of $(1\ |\ 3\ |\ 9)(5,6,7,8)$. In total, we obtain that the number of peaks in $D$ plus the number of descents in $\phi_n(D)$ add up to~$n$. Since the number of descents of $D$ is one less than the number of peaks, the lemma follows.
      \end{proof}
      Define the {\sf lifting} $\Delta:\D_{n} \rightarrow \D_{n+1}$ as $\Delta(D) \in \D_{n+1}$ being obtained from $D \in \D_{n}$ by adding a north step at the beginning of $D$ and an east step at the end.
      \begin{examplecon}
        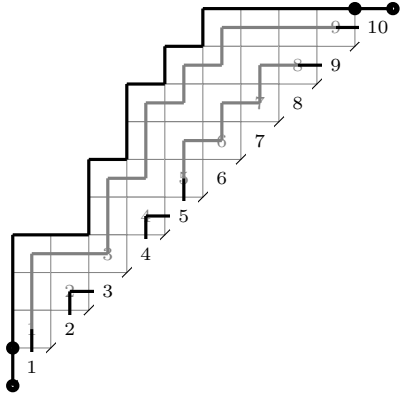
\begin{figure}
          \centering
          \begin{tikzpicture}[scale=1]
            \boxcollection{0.5}{(0,0)}{
            1/6/,1/7/,1/8/,
            2/5/,2/6/,2/7/,
            3/4/,3/5/,3/6/,
            4/4/,4/5/,
            5/3/,5/4/,
            6/3/,
            7/1/,7/2/,
            8/1/}
            \boxcollection{0.5}{(0,0)}{
            1/9/,2/8/,3/7/,4/6/,5/5/,6/4/,7/3/,8/2/,9/1/}
            \latticepath{0.5}{(0,-5.0)}{very thick}{
              0/4/black,2/0/black,0/2/black,1/0/black,0/2/black,1/0/black,0/1/black,1/0/black,0/1/black,5/0/black}
            \latticepath{0.5}{(0.25,-4.25)}{very thick}{
              0/2/grey,2/0/grey,0/2/grey,1/0/grey,0/2/grey,1/0/grey,0/1/grey,1/0/grey,0/1/grey,3/0/grey}
            \latticepath{0.5}{(2.25,-2.25)}{very thick}{
              0/1/grey,1/0/grey,0/1/grey,1/0/grey,0/1/grey,1/0/grey}
            \latticepath{0.5}{(0.0,-5.0)}{thin}{
              10/10/black}
            \content{0.5}{(0.5,-4.5)}{
               0/0/{\color{grey}$\scriptstyle 1$},
              -1/1/{\color{grey}$\scriptstyle 2$},
              -2/2/{\color{grey}$\scriptstyle 3$},
              -3/3/{\color{grey}$\scriptstyle 4$},
              -4/4/{\color{grey}$\scriptstyle 5$},
              -5/5/{\color{grey}$\scriptstyle 6$},
              -6/6/{\color{grey}$\scriptstyle 7$},
              -7/7/{\color{grey}$\scriptstyle 8$},
              -8/8/{\color{grey}$\scriptstyle 9$}}
            \latticepath{0.5}{(0.25,-4.75)}{very thick}{
              0/1/black}
            \latticepath{0.5}{(4.25,-0.25)}{very thick}{
              1/0/black}
            \latticepath{0.5}{(0.75,-4.25)}{very thick}{
              0/1/black,1/0/black}
            \latticepath{0.5}{(1.75,-3.25)}{very thick}{
              0/1/black,1/0/black,0/1/black}
            \latticepath{0.5}{(3.75,-0.75)}{very thick}{
              1/0/black}
            \content{0.5}{(0.5,-5.0)}{
              0/0/\colorbox{white}{$\scriptstyle 1$},-1/1/\colorbox{white}{$\scriptstyle 2$},-2/2/\colorbox{white}{$\scriptstyle 3$},-3/3/\colorbox{white}{$\scriptstyle 4$},-4/4/\colorbox{white}{$\scriptstyle 5$},-5/5/\colorbox{white}{$\scriptstyle 6$},-6/6/\colorbox{white}{$\scriptstyle 7$},-7/7/\colorbox{white}{$\scriptstyle 8$},-8/8/\colorbox{white}{$\scriptstyle 9$},-9/9.1/\colorbox{white}{$\scriptstyle 10$}}
             \draw[line width=3pt] (0.0,-5.0) circle (1pt);
             \draw[line width=3pt] (0.0,-4.5) circle (1pt);
             \draw[line width=3pt] (4.5,0.0) circle (1pt);
             \draw[line width=3pt] (5.0,0.0) circle (1pt);
          \end{tikzpicture}
          \caption{The lifting $\Delta(D)$ of the Dyck path $D$ in Figure~\ref{fig:bijNNNCA}.}
          \label{fig:mapdelta}
        \end{figure}
        The lifting $\Delta(D)$ of the Dyck path $D$ in Figure~\ref{fig:bijNNNCA} is shown in Figure~\ref{fig:mapdelta}.
        It is sent by the map $\phi_{10}$ to the noncrossing partition
        $$\sigma^{\Delta} = (1\ |\ 10)(2,3)(4,5\ |\ 9) = [10,3,2,5,9,6,7,8,4,1].$$
      \end{examplecon}

      \begin{lemma}\label{le:shift}
        Let $D \in \Dn$, $\sigma := \phi_n(D)$, $D^{\Delta} := \Delta(D)$, and $\sigma^{\Delta} := \phi_{n+1}(D^{\Delta})$. Then
        $$\maj(D^{\Delta}) + \maj(\sigma^{\Delta}) + \imaj(\sigma^{\Delta}) = \maj(D)+\maj(\sigma)+\imaj(\sigma)+2n.$$
      \end{lemma}
      \begin{proof}
        We first have to understand how the cycles of $\sigma^{\Delta}$ can be computed from the cycles of~$\sigma$.
        This is to say that we need to understand how the lifting $\Delta$ changes the structure of the cycles.
        Observe that via the lifting $\Delta$ and the shelling~$\phi$, a singleton~$i$ of~$\sigma$ correspond exactly to the pair $i,i+1$ of consecutive integers in a block of $\sigma^{\Delta}$, while a pair $j,j+1$ of consecutive integers in a block of~$\sigma$ exactly correspond to the singleton $j+1$ in $\sigma^{\Delta}$.
        This follows from the description of singletons and consecutive integers in blocks given in the proof of Lemma~\ref{le:descentvalleys}.
        For example, the singletons $2$ and $4$ in the permutation~$\sigma$ in the above example correspond to the consecutive integers $23$ and $45$ within cycles of $\sigma^{\Delta}$.
        On the other hand, the consecutive integers $56$, $67$, and $78$ correspond to the singletons $6$, $7$, and $8$ inside $\sigma^{\Delta}$.
        In other words, singletons of~$\sigma$ correspond exactly to the nonmaximal integers in the blocks of $\sigma^{\Delta}$, and nonminimal integers in the blocks of~$\sigma$ correspond exactly to the singletons in $\sigma^{\Delta}$.
        Using this observation, we now have one block in a cycle of $\sigma^{\Delta}$ that is formed by the integer $n+1$ together with all its preceding singletons in~$\sigma$, and that all other blocks in cycles of $\sigma^{\Delta}$ are given by the minimal elements in blocks of~$\sigma$ together with their preceding singletons in~$\sigma$.

        It is now left to understand which blocks are joined together to a cycle in $\sigma^{\Delta}$.
        This is obtained by observing that the block containing the smallest integer $i_1$ in a cycle $(i_1,\ldots,i_k)$ of~$\sigma$ is joined to the block containing the integer $i_k+1$.
        In the previous example, we had that $\sigma = (1\ |\ 3\ |\ 9)(5,6,7,8)$. Together with the singletons, we obtain that the blocks in the cycles of $\sigma^{\Delta}$ are given by
        $$1 \quad|\quad 2,3 \quad|\quad 9 \quad|\quad 4,5 \quad|\quad 10.$$
        Moreover, the block containing $1$ is joined to the block containing $10 = 9+1$, and the block containing $5$ is joined to the block containing $9 = 8+1$. We thus have in total that $\sigma^{\Delta}=(1\ |\ 10)(2,3)(4,5\ |\ 9)$, as expected.

        To proceed with the proof, we need the observation (which follows immediately from the previous considerations) that the minimal (resp. maximal) blocks in cycles of~$\sigma$ exactly correspond to the nonmaximal (resp. nonminimal) blocks in cycles of $\sigma^{\Delta}$.
        Together with the descriptions of the descent and inverse descent set, we now have that
        $$\Des(\sigma^{\Delta}) = \iDes(\sigma) \sqcup \{n\}, \quad \iDes(\sigma^{\Delta}) = \{i+1 : i \in \Des(\sigma)\} \sqcup \{1\},$$
        where we use the symbol $\sqcup$ to denote the disjoint union.
        Therefore,
        $$\maj(\sigma^{\Delta}) = \imaj(\sigma) + n, \quad \imaj(\sigma^{\Delta}) = \maj(\sigma)+\des(\sigma)+1.$$
        On the other hand,
        $$\maj(D^{\Delta}) = \maj(D) + \des(D).$$
        The statement follows with Lemma~\ref{le:descentvalleys}.
      \end{proof}
      \begin{examplecon}
        In our ongoing example, we have already seen that
        $$\maj(D) + \maj(\sigma) + \imaj(\sigma) + 2n = 72 + 18 = 90 = 9 \cdot 10.$$
        On the other hand, we have
        $$\Des(D^{\Delta}) = \{6,9,12,14\}, \quad \Des(\sigma^{\Delta}) = \{1,2,5,8,9\}, \quad \iDes(\sigma^{\Delta}) = \{1,2,4,8,9\}.$$
        This gives
        \begin{align*}
          \maj(D^{\Delta}) &= 6 + 9 + 12 + 14 = 41, \\
          \maj(\sigma^{\Delta}) + \imaj(\sigma^{\Delta}) &= (1+2+5+8+9) + (1+2+4+8+9) = 25 + 24, \\
          \maj(D^{\Delta}) + \maj(\sigma^{\Delta}) + \imaj(\sigma^{\Delta}) &= 41 + 25 + 24 = 90 = 9 \cdot 10.
        \end{align*}
      \end{examplecon}
      \begin{proof}[Proof of Proposition~\ref{prop:majimajA} for~$\phi_n$]
        We prove the proposition by induction on the length of the Dyck path. Let $D \in \Dn, D' \in \D_{n'}$ with $n,n' \geq 1$. By $DD'$, we denote the concatenation of $D$ and $D'$ in $\D_{n+n'}$. The proof consists of two parts,
        \begin{itemize}
            \item[(i)] if the theorem holds for $D$ and $D'$ then it holds for $DD'$, and
            \item[(ii)] if the theorem holds for $D$ then it holds for $\Delta(D)$.
        \end{itemize}
         As the case $n = 1$ is obvious, the theorem then follows.
        \begin{itemize}
          \item[(i)] set $\sigma:= \phi_n(D), \sigma' := \phi_{n'}(D')$ and $\tau := \phi_{n+n'}(DD')$.
          We then have
          \begin{align*}
            \maj(DD') &= \sum_{i \in \Des(DD')}i = \sum_{i\in\Des(D)}i \quad+\quad \sum_{i\in\Des(D')}(2n+i) \quad+\quad 2n\\
                      &= \maj(D)+\maj(D')+2n(\des(D')+1)
          \end{align*}
          where the summand $(2n+i)$ in the second sum comes from the additional length~$2n$ of $D$, and where the last~$2n$ is the additional descent between the end of $D$ and the beginning of $D'$.
          Moreover,
          \begin{align*}
            \maj(\tau) &= \sum_{i \in \Des(\tau)}i = \sum_{i\in\Des(\sigma)}i \quad+\quad \sum_{i\in\Des(\sigma')}(n+i)\\
                      &= \maj(\sigma)+\maj(\sigma')+n\des(\sigma')
          \end{align*}
          where we used the property of~$\phi$ that the one-line notation of~$\tau$ is obtained as the concatenation of the one-line notations of~$\sigma$ and of~$\sigma'$, with the one-line notation of~$\sigma'$ shifted by~$n$.
          In symbols,
          $$\tau_i = \begin{cases} \sigma_i &\text{ if } 1 \leq i \leq n, \\ \sigma'_{i-n} + n &\text{ if } n < i \leq n+n' \end{cases}.$$
          The same observation gives
          \begin{align*}
            \imaj(\tau) &= \sum_{i \in \iDes(\tau)}i = \sum_{i\in\iDes(\sigma)}i + \sum_{i\in\iDes(\sigma')}(n+i)\\
                      &= \imaj(\sigma)+\imaj(\sigma')+n\ides(\sigma').
          \end{align*}
          In total, we thus have
          \begin{align*}
            \maj(DD') + \maj(\tau) + \imaj(\tau)   &= \maj(D) + \maj(D') + 2n(\des(D')+1)\\
                                                  &+ \maj(\sigma)+\maj(\sigma') + n \des(\sigma')\\
                                                  &+ \imaj(\sigma) + \imaj(\sigma') + n \ides(\sigma').
          \end{align*}
          Lemma~\ref{le:NNdes} now gives $\des(\sigma') = \ides(\sigma')$ and $n\des(\sigma') + n\ides(\sigma') = 2n\des(\sigma')$, and Lemma~\ref{le:descentvalleys} then gives $2n\des(\sigma') + 2n(\des(D')+1) = 2nn'$.
          Putting this together finally gives
          \begin{align*}
            \maj(DD') & + \maj(\tau) + \imaj(\tau) = \\ & \maj(D)+\maj(\sigma)+\imaj(\sigma)+\maj(D')+\maj(\sigma') +\imaj(\sigma')+2nn'.
          \end{align*}
          By induction, this reduces to $n(n-1) + n'(n'-1) + 2nn' = (n+n')(n+n'-1)$.
          \item[(ii)] this is an immediate consequence of Lemma~\ref{le:shift}.
        \end{itemize}
      \end{proof}

    \subsection{Dyck paths and $231$-avoiding permutations}
    \label{sec:DyckCox}
    
    Bijections between $3$-pattern-avoiding permutations and Dyck paths are very well studied, see for example~\cite{BK2001, Cal2007, Kra2001, Rei2003, Stu2008}.
    Before proving Theorem~\ref{th:bijectionsA} for $231$-avoiding permutations, observe that Theorem~\ref{th:majimajNCCox} is in this case a direct consequence of \cite[Corollary~3.12]{Stu2008}.
    The advantage of studying an alternative bijection is that it also keeps track of the area statistic on Dyck paths, and that it will guide us later to obtain as well the analogous result for the group of signed permutations. The provided alternative bijection will be described using several bijections studied in~\cite{Stu2008}, see below.

    Let $D \in \Dn$ be a Dyck path of length~$2n$. Label every box $b_{ij}$ contributing to the area of $D$ by $s_{n-1-i}$ as shown in Figure~\ref{fig:dyck5Cox}. The bijection between Dyck paths and $231$-avoiding permutations is then defined by mapping $D$ to $\psi_n(D) := \prod s$, where the product ranges over all simple transpositions in the labeled boxes in the order as indicated in the figure. To see that $\psi_n(D)$ is in fact in $\Cox_n$, we use a description of $231$-avoiding permutations in \cite[Exercise~2.2.1.4--5]{Knu1973}, see also \cite[Example~2.3]{Rea2007}. A permutation $\sigma \in \Sn$ is $231$-avoiding if and only if it has a reduced expression of the form $X_1 X_2 \ldots X_{n-1}$ where each $X_\ell$ is a (possibly empty) subword of $s_{n-1}\cdots s_1$ and where moreover all simple transpositions in $X_{\ell+1}$ are also contained in $X_\ell$.
    \begin{example}\label{ex:dyck5Cox}
      The Dyck path shown in Figure~\ref{fig:dyck5Cox} is mapped to the $231$-avoiding permutation
      $$s_5s_4s_3s_2s_1 | s_5s_4s_2 | s_5 = [6,2,1,5,4,3] \in \Cox_6.$$
      \begin{figure}
        \centering
          \begin{tikzpicture}[scale=1]
            \boxcollection{0.5}{(0,0)}{
            1/1/,1/2/,1/3/,1/4/,1/5/,
            2/1/,2/2/,2/3/,2/4/,
            3/1/,3/2/,3/3/,
            4/1/,4/2/,
            5/1/}
            \latticepath{0.5}{(0,-3)}{very thick}{
              0/4/black,3/0/black,0/2/black,3/0/black}
            \latticepath{0.5}{(0,-3)}{very thin}{
              6/6/black}
            \latticepath{0.5}{(0,-2.5)}{very thin}{
              5/5/grey}
            \latticepath{0.5}{(0,-2)}{very thin}{
              4/4/grey}
            \latticepath{0.5}{(0,-1.5)}{very thin}{
              3/3/grey}
            \latticepath{0.5}{(0,-1)}{very thin}{
              2/2/grey}
            \latticepath{0.5}{(0,-0.5)}{very thin}{
              1/1/grey}
            \content{0.5}{(0.5,-2.5)}{
              0/0/$\color{black}s_5$,-1/1/$\color{black}s_4$,-2/2/$\color{black}s_3$,-3/3/$\color{black}s_2$,-4/4/$\color{black}s_1$,
              -1/0/$\color{black}s_5$,-2/1/$\color{black}s_4$,-4/3/$\color{black}s_2$,
              -2/0/$\color{black}s_5$}
            \content{0.5}{(-0.5,-1.5)}{
              0/0/$\scriptstyle 5$,1/0/$\scriptstyle 4$,2/0/$\scriptstyle 3$,3/0/$\scriptstyle 2$,4/0/$\scriptstyle 1$}
            \draw [->] (-0.5,-3.0) -- (-0.2,-2.7);
            \draw [->] (-0.5,-2.5) -- (-0.2,-2.2);
            \draw [->] (-0.5,-2.0) -- (-0.2,-1.7);
            \draw [->] (-0.5,-1.5) -- (-0.2,-1.2);
            \draw [->] (-0.5,-1.0) -- (-0.2,-0.7);
          \end{tikzpicture}
        \caption{A Dyck path of length $6$ with boxes labeled by simple transpositions.}
        \label{fig:dyck5Cox}
      \end{figure}
    \end{example}
    In the case of $231$-avoiding permutation, Theorem~\ref{th:bijectionsA} for $\psi_n:\Dn\tilde\longrightarrow\NC_n$ is a direct consequence of the following proposition, where $\Set_X(D)$ and $\Set_Y(D)$ denote the collection of~$x$- and~$y$-coordinates of the descents of $D$.
    \begin{proposition}\label{pr:importantbij}
      Let $D$ be a Dyck path and let $\sigma = \psi_n(D)$. Then $\area(D) = \inv(\sigma)$ and furthermore,
      \begin{align*}
        \Des(\sigma) &= [n-1] \setminus \{n-i : i \in \Set_X(D)\}, \\
        \iDes(\sigma) &= [n-1] \setminus \{n-i : i \in \Set_Y(D)\}.
      \end{align*}
      In particular, we obtain
      $$\maj(D) + \maj(\sigma) + \imaj(\sigma) = n(n-1).$$
    \end{proposition}
    \begin{proof}
      The fact that $\area(D) = \inv(\sigma)$ follows directly from the construction, as the inversion number $\inv(\sigma)$ is given by the minimal number of simple transpositions needed to write~$\sigma$.
      To obtain the statements about the descent and the inverse descent set of~$\sigma$, we follow several bijections described in~\cite{Stu2008} together with a bijection $f_n$ described by J.~Bandlow and K.~Killpatrick in~\cite[Lemma~1]{BK2001}.
      Observe that~$\psi_n$ and $f_n$ only differ by interchanging~$s_i$ and~$s_{n-i}$ in the definition (be aware of the different convention in~\cite{BK2001} where permutations are multiplied from left to right). Since $w_\circ s_i w_\circ = s_{n-i}$ for the involution $w_\circ := [n,\ldots,2,1]$, we thus have $\psi_n = \alpha_n \circ f_n$ where $\alpha_n$ is the involution on permutations that replaces $\sigma = [\sigma_1,\ldots,\sigma_n]$ by $\alpha_n(\sigma) = w_\circ \sigma w_\circ = [n+1-\sigma_n,\ldots,n+1-\sigma_1]$.
      At the end of~\cite[Section~5]{Stu2008} we have seen how to describe $f_n$ in terms of three other bijections, $f_n = {\bf i} \circ \beta_n \circ \gamma_n$, where
      \begin{itemize}
        \item ${\bf i}$ sends~$\sigma$ to $\sigma^{-1}$ (and thus interchanging the descent and the inverse descent sets),
        \item $\beta_n$ sends $\sigma \in \Cox_n$ to the unique Dyck path $D = \beta_n(\sigma)$ for which $\Des(\sigma) = \Set_X(D)$ and $\iDes(\sigma) = \Set_Y(D)$, and
        \item $\gamma_n$ sends a Dyck path $D$ to the unique Dyck paths $D' = \gamma(D)$ for which $\Set_X(D') = [n-1] \setminus \Set_Y(D)$ and $\Set_Y(D') = [n-1] \setminus \Set_X(D)$.
      \end{itemize}
      Thus, $\psi_n = \alpha_n \circ {\bf i} \circ \beta_n \circ \gamma_n$, and we obtain the sequence of equalities
      \begin{align*}
        \Des(\sigma)  &= \big\{ n-i : i \in \Des(\alpha(\sigma)) \big\} \\
                      &= \big\{ n-i : i \in \iDes({\bf i} \circ \alpha(\sigma)) \big\} \\
                      &= \big\{ n-i : i \in \Set_Y(\beta\circ{\bf i}\circ\alpha(\sigma)) \big\}
                      \\
                      &= \big\{ n-i : i \in [n-1] \setminus \Set_X(D) \big\}
                      \\
                      &= [n-1] \setminus \big\{ n-i : i \in \Set_X(D) \big\},
                      \\
      \end{align*}
      as desired. The sequence of equalities is similar for $\iDes(\sigma)$. Observe here, that the first equality is a direct property of the map $\alpha_n$.

      To finally see that $\maj(D) + \maj(\sigma) + \imaj(\sigma) = n(n-1)$, observe that $\maj(D)$ is equal to the sum of the entries in $\Set_X(D)$ plus the sum of the entires in $\Set_Y(D)$, and thus,
      $$\maj(D) + \maj(\sigma) + \imaj(\sigma) = 2\sum_{i \in [n-1]}i = n(n-1).$$
    \end{proof}
    We can as well keep track of another statistic on Dyck paths, namely the length of the last descent. We will use this fact along the way when describing generalizations of the constructions for signed permutations in the next section.
    \begin{proposition}\label{pr:lengthoflastdescent}
      Let $D$ be a Dyck path of length~$n$ and let $k$ be the number of east steps after the last north step. Then $\sigma_k = 1$ for $\sigma = \psi_n(D)$, and furthermore, $\{1,\ldots,k-1 \} \subseteq \Des(\sigma)$.
    \end{proposition}
    \begin{proof}
      Let $X_1 X_2 \cdots X_k$ be the initial segment of $\sigma = X_1 X_2 \cdots X_{n-1}$, with $X_k$ possibly empty. Then, by construction, the last simple transposition in $X_i$ is $s_i$ for $i < k$ and $s_k$ is \emph{not} contained in $X_k$. Therefore, $k$ is mapped by~$\sigma$ to $1$. As~$\sigma$ is $231$-avoiding, we also have $\{1,\ldots,k-1 \} \subseteq \Des(\sigma)$.
    \end{proof}

  \section{Definitions and proofs for the group of signed permutations} \label{sec:signedpermutations}

    In this section, we generalize the constructions in Section~\ref{sec:permutations} to signed permutations.

    \subsection{Dyck paths for signed permutations}

      The following definition of Dyck paths for signed permutations is motivated by the connection between Dyck paths and order ideals in the root poset of type $A$ as defined e.g. in \cite[Section~5.1.2]{Arm2006}. A \defn{type $B_n$ Dyck path} is a lattice path of~$2n$ steps, starting at $(0,0)$, consisting of north and east steps, and which never goes below the diagonal $x=y$. We denote the set of all type $B_n$ Dyck paths by $\D_{B_n}$. The \defn{area statistic} $\area(D)$ is defined to be the number of boxes $b_{ij}$ (as defined in~\eqref{eq:boxes} above) which lie below $D$ and for which $1 \leq i < j < 2n-i$. An example of a Dyck path in $\D_{B_6}$ is shown in Figure~\ref{fig:atypebcatpath}; all boxes contributing to its area are marked with little circles.

      \medskip

      It is easy to check, see \cite{FH1985}, that the area generating function for $\Dn$, $\Cat_n(q) = \sum_{D \in \Dn}q^{\area(D)}$, satisfies a recurrence given by
      $$\Cat_{n}(q) = \sum_{\ell=0}^{n-1} q^\ell \Cat_\ell(q) \Cat_{n-1-\ell}(q).$$
      For $\Cat_{B_n}(q) = \sum_{D \in \D_{B_n}}q^{\area(D)}$, an analogous statement holds.
      \begin{theorem}\label{th:recurrence}
        The $q$-Catalan numbers $\Cat_{B_n}(q)$ satisfy the recurrence relation
        $$\Cat_{B_n}(q) = \Cat_n(q) + \sum_{\ell=0}^{n-1}{q^{2\ell+1} \Cat_{B_\ell}(q) \Cat_{n-\ell}(q)}, \qquad \Cat_{B_0}(q) = 1.$$
      \end{theorem}
      \begin{proof}
        $D \in \D_{B_n}$ has either exactly~$n$ north and~$n$ east steps, which means it lies in $\Dn$, or there exists a last point $(\ell,\ell+1)$ where the path touches the diagonal $x+1=y$ and stays strictly above afterwards. Now, we have an initial Dyck path in $\D_{\ell+1}$, except that the last step is a north step instead of an east step, see Figure~\ref{fig:atypebcatpath} for an example. After this north step, a Dyck path in $\D_{B_{n-\ell-1}}$ starts. This gives
        \begin{align}
          \aCat_{B_n}(q)  &= \aCat_n(q) + \sum_{\ell=0}^{n-1}{\aCat_{\ell+1}(q) q^{2(n-\ell)-1}\aCat_{B_{n-\ell-1}}(q)} \label{eq:Bnrecurrence} \\
                          &= \aCat_n(q) + \sum_{\ell=0}^{n-1}{q^{2\ell+1} \aCat_{B_\ell}(q) \aCat_{n-\ell}(q)} \nonumber
        \end{align}
        \begin{figure}
          \centering
          \begin{tikzpicture}[scale=1]
            \boxcollection{0.5}{(0,0)}{
            1/4/,
            2/4/,2/5/,
            3/3/,3/4/,3/5/,
            4/3/,4/4/,
            5/3/,
            6/2/}
            \content{0.5}{(0,0)}{
            1/4/$\circ$,1/5/$\circ$,
            2/4/$\circ$,2/5/$\circ$,2/6/$\circ$,
            3/3/$\circ$,3/4/$\circ$,3/5/$\circ$,
            4/3/$\circ$,4/4/$\circ$,
            5/3/$\circ$,
            6/2/$\circ$}
             \latticepath{0.5}{(0,-4)}{very thick}{
               0/1/black,1/0/black,0/2/black,1/0/black,0/3/black,1/0/black,0/2/black,1/0/black}
             \latticepath{0.5}{(0,-4)}{very thin}{
               6/6/black,-2/2/black}
             \latticepath{0.5}{(3,-1)}{very thin}{
               -1/0/black}
             \latticepath{0.5}{(1,-2.5)}{very thin}{
               3.5/3.5/black}
             \latticepath{0.5}{(1,-2.0)}{very thin}{
               3/3/black}
             \draw[line width=3pt] (1,-2) circle (1pt);
             \draw[line width=3pt] (1,-2.5) circle (1pt);
          \end{tikzpicture}
          \caption{A Dyck path in $\D_{B_6}$ with $\area$ $12$.}
          \label{fig:atypebcatpath}
        \end{figure}
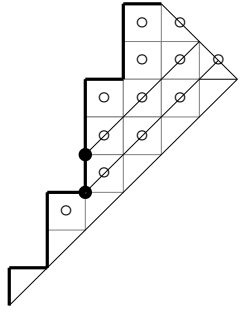
        The factor $q^{2(n-\ell)-1}$ in~\eqref{eq:Bnrecurrence} comes from the difference between the area of the given type $B_n$ Dyck path and the sum of the areas of its decomposition into an (almost) Dyck path in $\D_{\ell+1}$ and a Dyck path in $\D_{B_{n-\ell-1}}$.
      \end{proof}
    \begin{example}
      Figure~\ref{fig:atypebcatpath} shows a Dyck path in $\D_{B_6}$. It starts with a path which is almost a Dyck path in $\D_3$, ending with the north step between the two dots, followed by a Dyck path in $\D_{B_3}$, which starts at the second dot. Thus, $n = 6$ and $\ell = 2$ in this case. The area of the the given Dyck path is $12$, which is decomposed into the area of the almost Dyck path in $\D_3$, which is $1$, the area of the Dyck path in $\D_{B_3}$, which is $4$, and the additional $7 = 2(6-2)-1$ boxes on the two diagonals starting at the two dots.
    \end{example}
    Let $(a;q)_n$ be shorthand for the $q$\emph{-shifted factorial} $\prod_{\ell=0}^{n-1}(1-aq^\ell)$.
    \begin{corollary}\label{corgenfunctionB}
        $\aCat_{B_n}(q)$ satisfies the generating function identity
        $$\sum_{n \geq 0}{\frac{x^n q^{-n(n-1)}(1-qx)}{(-x;q^{-1})_{2n+1}}\aCat_{B_n}(q)} = 1.$$
    \end{corollary}
    \begin{proof}
      The recurrence in Theorem~\ref{th:recurrence} can be written as
      $$
        (1+q^{2n+1})\aCat_{B_n}(q) -  \sum_{\ell=0}^{n}{q^{2\ell+1} \aCat_{B_\ell}(q) \aCat_{n-\ell}(q)} = \aCat_n(q).
      $$
      Multiplying both sides of the equality by $\frac{{x^n q^{-n(n-1)}}}{{(-x;q^{-1})_{2n+1}}}$ and summing over all~$n$ yields the proposed recurrence relation.
      Here, we used that in~\cite[Theorem~5]{FH1985} with $\aCat_n(q) = C_n(q;q^2,q^{-2})$, it is shown that $\aCat_{n}(q)$ satisfies the generating function identity
      $$x = \sum_{n \geq 1} \frac{(x+1)x^nq^{-n(n-1)}}{(-x;q^{-1})_{2n+1}}\aCat_n(q),$$
      and thus
      $$\sum_{n\geq 0}\frac{x^nq^{-n(n-1)}}{(-x;q^{-1})_{2n+1}}\aCat_n(q) = 1.$$
    \end{proof}
    After discussing some basic properties of the area generating function of type $B_n$ Dyck paths, we now turn to a version of the major index for Dyck paths of type $B_n$.
    In this case, we consider~$2n$ as well a descent of $D$ if the $(2n)$-th step of $D$ (which is the step at its loose end) is an east step.
    The reason for this convention will become clear in the subsequent constructions.
    In other words, the \defn{flag descent set} $\fDes(D)$ is the usual descent set together with this possible descent at position~$2n$, and $\fdes(D) = |\fDes(D)|$.
    We moreover define the \defn{flag major index} of $D$ as
    $$\fmaj(D) := 2 \cdot \bigg( \neg(D) + \sum_{i \in \fDes(D)}{(2n-i)} \bigg),$$
    where $\neg(D)$ equals the number of east steps in $D$.
    Observe that~$i$ is replaced by $2n-i$, as $D \in \D_{B_n}$ should be considered to \lq\lq start\rq\rq\ at the loose end.
    In particular, a possible descent at this loose end does not contribute to the flag major index.

    For example, the descent set of the Dyck path $D = NENNENNNENNE \in \D_{B_6}$ shown in Figure~\ref{fig:atypebcatpath} is given by $\Des(D) = \{ 2,5,9,12 \}$, and the flag major index thus equals
    $$\fmaj(D) = 2( 4 + (12-2) + (12-5) + (12-9) + (12-12) ) = 48.$$
    We chose the term \lq\lq flag major index\rq\rq\ in analogy with the flag major index for signed permutations, which we will discuss in Section~\ref{sec:DyckReverseNC}.

    One can define the flag major index alternatively using the peaks of $D$ rather than the valleys. The reason, we used the valleys in the definition is the similarity with the definition of the flag major index of permutations in type $B_n$. Let $\Asc(D)$ is the number of indices~$i$ for which the~$i$-th step of $D$ is a north step and the $(i+1)$-st step is an east step. We then have the following lemma.
    In the subsequent constructions and proofs, we will often change back and forth between the definition of the flag major index, and the below description in terms of ascents.
    \begin{lemma}\label{le:peaks}
      Let $D$ be a Dyck path of type $B_n$. Then
      $$\fmaj(D) = 2\sum_{i\in\Asc(D)}(2n-i).$$
    \end{lemma}
    \begin{proof}
      There is a one-to-one correspondence between peaks (or ascents) of $D$ and its valleys (or descents).
      Even the last peak has always a following valley, since we consider~$2n$ to be a valley if it is an east step.
      Moreover, let $e_i$ be the number of east steps between the~$i$-th peak and the following~$i$-th valley. Then
      $$\sum_{i\in\Asc(D)}(2n-i) = \sum_{i\in\Des(D)}(2n-i) \quad+\quad \sum e_i = \sum_{i\in\Des(D)}(2n-i) \quad+\quad \neg(D).$$
      In the last equality, we used that every east step in $D$ contributes to exactly one of the $e_i$'s.
    \end{proof}

    \begin{proposition}\label{pr:typeBdCat}
      The generating function for the flag major index on Dyck paths in $\D_{B_n}$ is given by MacMahon's $q$-Catalan numbers for $B_n$,
      $$\sum_{D \in \D_{B_n}} q^{\fmaj(D)} = \begin{bmatrix} 2n \\ n \end{bmatrix}_{q^2}.$$
    \end{proposition}
    \begin{figure}
      \centering
          \begin{tikzpicture}[scale=1]
            \boxcollection{0.5}{(0,0)}{
            1/1/,1/2/,1/3/,1/4/,1/5/,1/6/,
            2/1/,2/2/,2/3/,2/4/,2/5/,2/6/,
            3/1/,3/2/,3/3/,3/4/,3/5/,3/6/,
            4/1/,4/2/,4/3/,4/4/,4/5/,4/6/,
            5/1/,5/2/,5/3/,5/4/,5/5/,5/6/,
            6/1/,6/2/,6/3/,6/4/,6/5/,6/6/}
            \latticepath{0.5}{(0,-3)}{thick}{
              0/1/black,2/0/black,0/1/black,2/0/black,0/2/black,1/0/black,0/2/black,1/0/black}
            \latticepath{0.5}{(0,-3)}{very thin}{
              6/6/grey}
            \latticepath{0.5}{(0.5,-3)}{very thin}{
              5/5/grey}
            \latticepath{0.5}{(1,-3)}{very thin}{
              4/4/grey}
            \latticepath{0.5}{(2,-2)}{line width=2.5pt}{
              -1/0/grey}
            \latticepath{0.5}{(1,-2.5)}{line width=2.5pt}{
              -1/0/grey}
             \draw[line width=3pt] (0,-2.5) circle (1pt);
             \draw[line width=3pt] (1,-2) circle (1pt);
             \draw[line width=3pt] (2,-1) circle (1pt);
             \draw[line width=3pt] (2.5,0) circle (1pt);
          \end{tikzpicture}
      \caption{The lattice path from $(0,0)$ to $(6,6)$ which is mapped to $D\in \D_{B_6}$ in Figure~\ref{fig:atypebcatpath}.}
      \label{fig:atypeblatticepath}
    \end{figure}
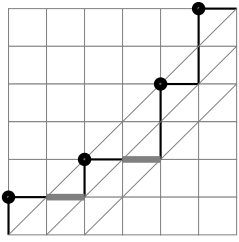
    \begin{proof}
      It is well-known that the major index generating function on lattice paths consisting of~$n$~north and~$n$~east steps without any further restrictions is given by
      $$\sum q^{\maj(L)} = \begin{bmatrix} 2n \\ n \end{bmatrix}_q,$$
      see for example~\cite{And1976} or~\cite{Mah1960}.
      By symmetry, we can assume that $\maj(L) = \sum_{i \in \Asc(L)} (2n-i)$. E.g., the ascents of the lattice path in Figure~\ref{fig:atypeblatticepath} are dotted, and we get $\Asc(L) = \{1,4,8,11\}$. Thus, the major index of $L$ is given by
      $$\maj(L) = (12-1)+(12-4)+(12-8)+(12-11) = 11+8+4+1 = 24.$$

      Define a bijection between lattice paths from $(0,0)$ to $(n,n)$ to Dyck paths in $\D_{B_n}$ by replacing the first east step from level~$i$ to level $i-1$ by a north step for all $i < 0$ for which such an east step exists. For example, the lattice path shown in Figure~\ref{fig:atypeblatticepath} is mapped to the Dyck path shown in Figure~\ref{fig:atypebcatpath}, the east steps which are replaced by north steps are drawn in {\color{grey}bold grey}. To see that this is indeed a bijection, observe that the inverse map is given by replacing the last north step from level~$i$ to level $i+1$ by an east step for all $0 \leq i < n-\neg(D)$ in a type $B_n$ Dyck path $D$.
      This transformation does not change the ascent set $\Asc$ and the sum $\sum_{i \in \Asc(L)}(2n-i)$ becomes $\sum_{i \in \Asc(D)}(2n-i)$ where $D$ is the type $B_n$ Dyck path obtained from $L$ via this transformation.
      Using Lemma~\ref{le:peaks}, we thus get $\fmaj(D) = 2 \maj(L)$ and the proposition follows.
    \end{proof}

  \subsection{Dyck paths and reverse noncrossing partitions}\label{sec:DyckReverseNC}

    A \defn{set partition $P$ of type $B_n$} is a set partition of $\{ \pm 1,\ldots,\pm n\}$ for which $B$ is a block of $P$ if and only if $-B$ is.
    Usually, one assumes as well that there is at most one central block $B$, i.e., a block $B$ for which $B = -B$. This comes from the fact that this is needed in order to obtain a combinatorial model for the intersection lattice of a certain type $B$ hyperplane arrangement, see e.g.~\cite{RS2010} for further details.
    We drop this condition here since we will consider \emph{reversed noncrossing partitions} of type $B_n$ which might contain multiple central blocks, see below.
    $P$ is noncrossing if it is noncrossing when $\{ \pm 1,\ldots,\pm n\}$ is ordered as
    \begin{align}
      -1<-2<\ldots<-n<1<2<\ldots<n. \label{eq:nonstandardordering}
    \end{align}
    We remark that the condition of being noncrossing implies that there is at most one central block. In particular, the definition of noncrossing partitions of type $B_n$ given here coincides with the usual definition.
    Following~\cite{RS2010}, the above ordering is, up to cyclic rotation, called \defn{crossing order}.
    One can naturally embed all set partitions of type $B_n$ into the set $B_n$ of signed permutations by sending a set partition $P$ to the permutation whose cycles are the blocks of $P$ in increasing order with respect to the above ordering.
    E.g., the noncrossing set partition
    $$\big\{ \{1\},\{-1\},\{ 2,3,-2,-3\},\{4,5\},\{-4,-5\},\{6,-6\} \big\}$$
    of type $B_6$ is sent to the signed permutation
    $$(2,3,-2,-3)(4,5)(-4,-5)(6,-6) = [1,3,-2,5,4,-6] \in B_6.$$
    We denote the set of noncrossing partitions of type $B_n$, considered inside $B_n$, by $\NC_{B_n}$.
    This definition can as well be found in~\cite{Rei2003} and in~\cite{Arm2006}.
    For our purposes, it is more convenient to work with the standard order
    \begin{align}
      -n<-(n-1)<\ldots<-1<1<2<\ldots<n. \label{eq:standardordering}
    \end{align}
    rather than with the crossing order.
    This order is called \defn{nesting order} in~\cite{RS2010}, and appears to be more natural when studying bijections between Dyck paths of type $B_n$ and noncrossing partitions of type $B_n$.
    We thus define the set $\rev(\NC_{B_n})$ of \defn{reverse noncrossing partitions} to be the set of all set partitions of type $B_n$ that are noncrossing with respect to the standard (or nesting) order.
    It is straightforward to check that
    $$\rev(\NC_{B_n}) = \big\{ \rev(\sigma) : \sigma \in \NC_{B_n} \big\},$$
    where $\rev$ is the involution on signed permutations which reverses all negative integers in the one-line notation. For example $\rev([2,-4,3,-1]) = [2,-1,3,-4]$.
    Set
    $$\Neg(\sigma) = \{ 1 \leq i \leq n : \sigma_i < 0\}$$
    to the the set of all positive indices~$i$ for which $\sigma_i$ is negative, and let $\neg(\sigma) = \big|\Neg(\sigma)\big|$.
    The \defn{flag descent set} for $\sigma = [\sigma_1,\ldots,\sigma_n] \in B_n$ is given by
    $$\fDes(\sigma) = \{ 0 \leq i < n : \sigma_i > \sigma_{i+1} \}.$$
    Here, we set $\sigma_0 = 0$, or equivalently, we consider $0$ to be a flag descent of~$\sigma$ if $1 \in \Neg(\sigma)$.
    We moreover set $\fdes(\sigma) = \big|\fDes(\sigma)\big|$ to be the number of flag descents of~$\sigma$.
    The \defn{flag major index} of~$\sigma$ is defined in \cite[Section~4]{ABR2001} as
    $$\fmaj(\sigma) := 2 \maj(\sigma) + \neg(\sigma),$$
    where $\maj(\sigma) = \sum_{i \in \fDes(\sigma)} i$.
    Observe that a possible flag descent at position $0$ does not contribute to the flag major index.
    Similarly, we again set $\iDes(\sigma) = \Des(\sigma^{-1})$, $\ifdes(\sigma) = \big|\ifDes(\sigma)\big|$, and $\ifmaj(\sigma) = \fmaj(\sigma^{-1})$.
    \begin{remark}
      As discussed in~\cite[Section~4]{ABR2001}, one can define the flag major index using either the standard or the crossing order. We used here the standard order since this approach seems to be better suited for connections to Dyck paths of type $B_n$. On the other hand, if one uses the crossing order instead, then the involution $\rev$ would yield a version of Corollary~\ref{cor:majimajB} with the flag major index for the crossing order (denoted below by $\fmaj^*(\sigma)$), and for noncrossing partitions rather than reversed noncrossing partitions,
      \begin{align*}
        \sum_{\sigma \in \NC_{B_n}} q^{\fmaj^*(\sigma) + \ifmaj^*(\sigma)} =  \begin{bmatrix} 2n \\ n \end{bmatrix}_{q^2}.
      \end{align*}

    \end{remark}
      \begin{figure}
        \centering
          \begin{tikzpicture}[scale=1]
            \boxcollection{0.5}{(0,0)}{
            0/2/$\circ$,
            1/2/$\circ$,1/3/$\circ$,
            2/1/$\circ$,2/2/$\circ$,2/3/$\circ$,
            3/1/$\circ$,3/2/$\circ$,
            4/-1/$\circ$,4/0/$\circ$,4/1/$\circ$,
            5/-1/$\circ$,5/0/$\circ$,
            6/-1/$\circ$}
            \content{0.5}{(0,0)}{
            0/3/$\circ$,
            1/4/$\circ$}
             \latticepath{0.5}{(-1,-3.5)}{very thick}{
               0/4/black,2/0/black,0/2/black,1/0/black,0/2/black,1/0/black}
             \latticepath{0.5}{(-0.75,-3.25)}{very thick}{
               0/3/black,2/0/black,0/2/black,1/0/black,0/2/black,1/0/black,0.15/0.2/black}
             \latticepath{0.5}{(-0.25,-2.75)}{very thick}{
               0/1/black,1/0/black}
             \latticepath{0.5}{( 0.75,-1.75)}{very thick}{
               0/1/black,1/0/black,0/2/black,1/0/black,0.15/0.2/black}
             \latticepath{0.5}{(-1,-3.5)}{very thin}{
               6/6/black,-2/2/black}
             \latticepath{0.5}{(2,-0.5)}{very thin}{
               -1/0/black}
            \content{0.5}{(0.5,-2.5)}{
              2/-2/\colorbox{white}{$\scriptstyle 6$},1/-1/\colorbox{white}{$\scriptstyle 5$},0/0/\colorbox{white}{$\scriptstyle 4$},-1/1/\colorbox{white}{$\scriptstyle 3$},-2/2/\colorbox{white}{$\scriptstyle 2$},-3/3/\colorbox{white}{$\scriptstyle 1$}}
          \end{tikzpicture}
        \caption{The bijection $\phi_{B_6}$ sending the shown Dyck path in $\D_{B_6}$ to the reverse noncrossing partition $(2,3,-2,-3)(4,5)(6,-6) \in \rev(\NC_{B_6})$.}
        \label{fig:NNNCB}
      \end{figure}
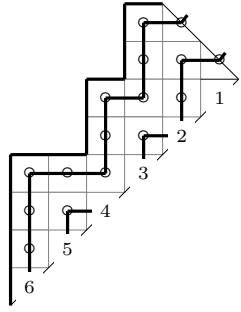
      To prove Theorem~\ref{th:bijectionsB}, we adapt the bijection $\phi_n : \Dn \longrightarrow \NC_n$ to signed permutations as follows. Write the numbers $1$ through~$n$ on the diagonal below the Dyck path $D \in \D_{B_n}$ as shown in Figure~\ref{fig:NNNCB}. Define $\phi_{B_n}(D)$ in the same way as for permutations with the additional rule that, if a \lq\lq shell\rq\rq\ ends at the \lq\lq right boundary\rq\rq, one adds the negatives of the elements to the cycle. This map defines a bijection between $\D_{B_n}$ and $\rev(\NC_{B_n})$.
      \begin{proposition}
        $\phi_{B_n}(D) : \D_{B_n} \longrightarrow \rev(\NC_{B_n})$ is well-defined and a bijection.
      \end{proposition}
      \begin{proof}
        Observe that a mirror reflection at the loose end of a Dyck path of type $B_n$ provides a centrally symmetric Dyck path of length $4n$, and that every centrally symmetric Dyck path of length $4n$ can be obtained this way. The map $\phi_{B_n}$ can now be considered as the original map $\phi_{2n}$ on the centrally symmetric Dyck path, where the integers $1$ through~$n$ are as well reflected to obtain the standard order (as we think of a type $B_n$ Dyck path to start from the loose end).
        As $\rev(\NC_{B_n})$ are the type $B_n$ set partitions that are noncrossing with respect to that ordering, we finally observe that the property of a Dyck path to be centrally symmetric exactly corresponds to the symmetry property of a set partition of type $B_n$ that $B$ is a block if and only if $-B$ is a block.
      \end{proof}

      For $\sigma \in B_n$, the \defn{type $B_n$ inversion number} $\inv_B(\sigma)$ is defined as the minimal $k$ for which there are $k$ simple transpositions in $\big\{ s_0 = (1,-1),s_1 = (1,2),\ldots, s_{n-1} = (n-1,n) \big\}$ whose product is~$\sigma$. It is well-known, see for example~\cite[Section~2]{ABR2001}, that
      $$\inv_B(\sigma) = \inv([\sigma_1,\sigma_2,\ldots,\sigma_n]) - \sum_{ i \in \Neg(\sigma)}\sigma_i,$$
      where $\inv$ is the usual inversion number in the one-line notation of~$\sigma$ computed with respect to the standard order.
      The following proposition thus proves Theorem~\ref{th:bijectionsB}\eqref{eq:areatoinvB} for $\phi_{B_n}:\D_{B_n} \tilde\longrightarrow\NC_{B_n}$.
      \begin{proposition}
        For $D \in \D_{B_n}$ and $\sigma = \phi_{B_n}(D)$, we have
        $$\area(D) = \inv([\sigma_1,\sigma_2,\ldots,\sigma_n]) - \sum_{i \in \Neg(\sigma)}\sigma_i.$$
      \end{proposition}
      \begin{proof}
        The proof follows exactly the same lines as the proof for permutations. First observe that if $D \in D_n \subseteq D_{B_n}$, then the statement reduces to Proposition~\ref{prop:area}. Now, assume that $D \in \D_{B_n} \setminus \D_n$ consists of a unique nontrivial cycle $\sigma = \phi_{B_n}(D) = (i_1,\ldots,i_k,-i_1,\ldots,-i_k)$. Then, the one-line notation of~$\sigma$ is obtained from the one-line notation of $\sigma' = (i_1,\ldots,i_k)$ by replacing~$i_1$ by~$-i_1$. Moreover, Let $D'$ be the preimage of~$\sigma'$ inside $\D_n \subseteq \D_{B_n}$.
        We then obtain that $\area(D) - \area(D') = 2i_1-1$. On the other hand, all inversions of~$\sigma'$ are as well inversions of~$\sigma$. Since all $1 \leq i < i_1$ appear to the left of $i_1$ in the one-line notation of~$\sigma'$,~$\sigma$ now has $i_1-1$ additional inversions.
        Finally,~$\sigma$ has exactly one negative entry, namely $i_1$.
        We thus get in total that
        $$\area(D) = \area(D') + 2i_1-1 = \inv(\sigma') + 2i_1-1 = \inv(\sigma) + i_1,$$
        as desired.
        As in the proof of Proposition~\ref{prop:area}, the general case is obtained by applying the same argument to each cycle of the form $(i_1,\ldots,i_k,-i_1,\ldots,-i_k)$ and to each pair of cycles of the form $(i_1,\ldots,i_k)(-i_1,\ldots,-i_k)$ individually.
      \end{proof}
      \begin{example}
        In Figure~\ref{fig:NNNCB}, the area of $D$ is $16$, and
        $$\sigma = \phi_{B_n}(D) = (2,3,-2,-3)(4,5)(6,-6) = [1,3,-2,5,4,-6].$$
        We therefore have
        $$\inv(\sigma) = \big|\{(1,3),(1,6),(2,3),(2,6),(3,6),(4,5),(4,6),(5,6)\}\big| = 8, \quad \Neg(\sigma) = \{3,6\}$$
        and thus $\area(D) = 8-(-2-6) = 16$, as expected.
      \end{example}

      The following proposition finally implies Theorem~\ref{th:bijectionsB}\eqref{eq:majtomajB} for $\phi_{B_n}:\D_{B_n} \tilde\longrightarrow\NC_{B_n}$.
      \begin{proposition}\label{propNNNCB}
        Let $D \in \D_{B_n}$. Then
        $$\fmaj(D) + \fmaj(\phi_{B_n}(D)) + \ifmaj(\phi_{B_n}(D)) = 2n^2.$$
      \end{proposition}
      \begin{proof}
        As in the proof of Proposition~\ref{prop:majimajA} for~$\phi_n$, we prove the statement by induction on the length of the Dyck path of type $B$.
        Let $D \in \D_n \subseteq \D_{B_n}$ and $D' \in \D_{B_{n'}}$ with $n \geq 1$ and $n' \geq 0$, and let $\tilde D$ be obtained from $D$ by replacing the last east step by a north step. Moreover, let $\tilde D D'$ denote the concatenation of $\tilde D$ and $D'$ in $\D_{B_{n+n'}}$.
        Following the decomposition of Dyck paths of type $B_n$ as discussed in (the proof of) Theorem~\ref{th:recurrence}, we have to prove two cases:
        \begin{itemize}
          \item[(i)] the theorem holds for $D$,
          \item[(ii)] if the theorem holds for $\tilde D$ and for $D'$, then it holds for $\tilde DD'$.
        \end{itemize}
        As indicated above, the decomposition in Theorem~\ref{th:recurrence} might leave $D'$ empty, i.e., $n' = 0$. We thus replace (ii) by
        \begin{itemize}
          \item[(iia)] the theorem holds for $\tilde D$,
          \item[(iib)] if the theorem holds for $D'$ with $n' \geq 1$, then it holds for $\tilde DD'$.
        \end{itemize}
        As it is again trivial to check that the statement holds for the two Dyck paths in $\D_{B_1}$, the theorem then follows.
        \begin{itemize}
          \item[(i)] Observe that $\phi_{B_n}(D) = \phi_n(D'')$, where $D'' = \rev(D)$ and $\rev : \D_n \longrightarrow \D_n$ interchanges north and east steps and reads the path backwards. The involution $\rev$ here accounts for the different conventions that the integer labellings on the diagonals for~$\phi_n$ start on the bottom left, while they start of the top right for $\phi_{B_n}$. We thus get
          \begin{align*}
            \fmaj(D) + \fmaj(\sigma) + \ifmaj(\sigma)
              &= 2\big( \maj(D'') + n + \maj(\sigma'') + \imaj(\sigma'') \big) \\
              &= 2( n(n-1) + n ) = 2n^2.
          \end{align*}
          Here, we set $\sigma'' = \phi_n(D'')$ and we used that $D$ and $D''$ have exactly~$n$ east steps, and that Proposition~\ref{prop:majimajA} yields
          $$\maj(D'') + \maj(\phi(D'')) + \imaj(\phi(D'')) = n(n-1).$$

          \item[(iia)] For more readability, we set
          $$\sigma = \phi_{B_n}(D), \quad \tilde\sigma = \phi_{B_n}(\tilde D), \quad \sigma' = \phi_{B_{n'}}(D'), \text{ and } \tau = \phi_{B_{n+n'}}(\tilde DD').$$
          First, observe that $\Neg(\sigma)$ is empty, and that $\tilde\sigma$ is obtained from~$\sigma$ by replacing $1$ by $-1$ in the one-line notation. In particular, this gives
          $\neg(\sigma) = 0$, and $\neg(\tilde\sigma) = 1$.

          We now have to distinguish two slightly different situations. Either $D$ ends in a north step followed by an east step (which results in an ascent at position $2n-1$), or in two east steps (where there is no ascent at position $2n-1$).

          In the first situation, we have $\sigma_1 = 1$ and thus $\tilde\sigma_1 = -1$. This gives
          \begin{align*}
            \Asc(\tilde D) \sqcup \{ 2n-1 \}  &= \Asc(D), \\
            \fDes(\tilde\sigma)  &= \fDes(\sigma) \sqcup \{ 0 \}, \\
            \ifDes(\tilde\sigma)  &= \ifDes(\sigma) \sqcup \{ 0 \},
          \end{align*}
          where we use the symbol $\sqcup$ to denote the disjoint union.
          Putting this together yields
          $$
            \fmaj(\tilde D) = \fmaj(D) - 2, \quad  \fmaj(\tilde\sigma)  = \fmaj(\sigma) + 1, \text{ and } \ifmaj(\tilde\sigma) = \ifmaj(\sigma) + 1.
          $$
          Using (i), the statement thus follows.

          In the second situation, we have that $\sigma_1 \neq 1$ and thus $\tilde\sigma_1 > 0$, while $(\tilde\sigma^{-1})_1 < 0$. This gives
          \begin{align*}
            \Asc(\tilde D)  &= \Asc(D), \\
            \fDes(\tilde\sigma)  &= \fDes(\sigma), \\
            \ifDes(\tilde\sigma) \sqcup \{ 0 \} &= \ifDes(\sigma) \sqcup \{ 1 \}.
          \end{align*}
          This now yields
          $$
            \fmaj(\tilde D) = \fmaj(D), \quad \fmaj(\tilde\sigma) = \fmaj(\sigma) + 1, \quad \ifmaj(\tilde\sigma) = \ifmaj(\sigma) + 1 - 2.
          $$
          Again using (i), the statement follows.

          \item[(iib)] We have
          \begin{align*}
            \fmaj(\tilde DD')  &= 2\sum_{i \in \Asc(\tilde DD')}(2(n+n')-i) \\
                        &= 2\sum_{i\in\Asc(D')}(2(n+n')-i-2n) \quad+\quad 2\sum_{i\in\Asc(\tilde D)}(2(n+n')-i) \\
                        &= \fmaj(D')+\fmaj(\tilde D)+4n'\fdes(\tilde D).
          \end{align*}
          Here, we used that $|\Asc(\tilde D)| = \fdes(\tilde D)$ and thus,
          $$\sum_{i\in\Asc(\tilde D)}(2(n+n')-i) = \sum_{i\in\Asc(\tilde D)}(2n-i) \quad+\quad 2n'\fdes(\tilde D).$$
          Next, it follows from the definition of the map $\phi_{B}$ that the one-line notation of~$\tau$ is obtained as the concatenation of the one-line notations of~$\sigma'$ and of $\tilde \sigma$, with the one-line notation of $\tilde \sigma$ shifted by $n'$. In symbols,
          $$\tau_i = \begin{cases} \sigma'_i &\text{ if } 1 \leq i \leq n', \\ \tilde\sigma_{i-n'} + n' &\text{ if } n' < i \leq n+n' \end{cases}.$$
          Therefore,
          \begin{align*}
            \fmaj(\tau) &= 2\maj(\tau)+\neg(\tau) = 2\sum_{i\in\fDes(\tau)}i \quad+\quad \neg(\tau) \\
                        &= 2\sum_{i\in\fDes(\sigma')}i \quad+\quad \neg(\sigma') \quad+\quad 2\sum_{i\in\fDes(\tilde\sigma)}(2n'+i) \quad+\quad \neg(\tilde\sigma) \\
                        &= \fmaj(\sigma')+\fmaj(\tilde\sigma)+2n'\fdes(\tilde\sigma).
          \end{align*}
          Here, we used that $0 \in \fDes(\tilde\sigma)$ if $\tilde\sigma_1 < 0$, and that $2n' \in \Des(\tau)$ if and only if $\tau_{2n'+1} = \tilde\sigma_1 < 0$.
          A similar consideration yields
          \begin{align*}
            \ifmaj(\tau) &= \ifmaj(\sigma')+\ifmaj(\tilde\sigma)+2n'\ides(\tilde\sigma).
          \end{align*}
          In total, we thus have
          \begin{align*}
            \fmaj(\tilde DD') + \fmaj(\tau) + \ifmaj(\tau) &= \fmaj(D')+\fmaj(\tilde D)+4n'\fdes(\tilde D) \\
                                                    &+ \fmaj(\sigma') + \fmaj(\tilde\sigma) + 2n'\fdes(\tilde\sigma) \\
                                                    &+ \ifmaj(\sigma') + \ifmaj(\tilde\sigma) + 2n'\ifdes(\tilde\sigma).
          \end{align*}
          Considering again the two cases for $\tilde D$ as above in (iia), we have
          \begin{align*}
            \fdes(\tilde D)       &= \fdes(D)-1 = \des(D), \\
            \fdes(\tilde\sigma)   &= \fdes(\sigma)+1 = \des(\sigma)+1 ,\\
            \ifdes(\tilde\sigma)  &= \ifdes(\sigma)+1 = \ides(\sigma)+1
          \end{align*}
          if $D$ has an ascent at position $2n-1$, and
          \begin{align*}
            \fdes(\tilde D)       &= \fdes(D) = \des(D)+1, \\
            \fdes(\tilde\sigma)   &= \fdes(\sigma) = \des(\sigma),\\
            \ifdes(\tilde\sigma)  &= \ifdes(\sigma)+1 = \des(\sigma)
          \end{align*}
          otherwise.
          In both cases, we can now use Lemma~\ref{le:NNdes} and Lemma~\ref{le:descentvalleys} for $D \in \D_n$, and obtain
          \begin{align*}
            4n'\fdes(\tilde D) + 2n'\fdes(\tilde\sigma) + 2n'\ifdes(\tilde\sigma)
              &= 4nn'.
          \end{align*}
          Thus, $\fmaj(\tilde DD') + \fmaj(\tau) + \ifmaj(\tau)$ equals
          $$
            \fmaj(\tilde D)+\fmaj(\tilde \sigma)+\ifmaj(\tilde \sigma)+\fmaj(D')+\fmaj(\sigma') +\ifmaj(\sigma')+4nn'.
          $$
          By induction, this reduces to $2n^2 + 2n'^2 + 4nn' = 2(n+n')^2$.
          \end{itemize}
    \end{proof}

  \subsection{Dyck paths and sortable elements}

    Following the general definition of Coxeter sortable elements, let $\Cox_{B_n}$ denote the set of signed permutations having a reduced expression of the form $X_1 X_2 \ldots X_n$ where each $X_\ell$ is a (possibly empty) subword of $s_{n-1}\cdots s_1 s_0$ and where moreover all simple transpositions in $X_{\ell+1}$ are also contained in $X_\ell$. As before, $s_0 = (1,-1) \in B_n$ is the signed permutation interchanging $1$ and $-1$, and $s_i = (i,i+1)$ for $1 \leq i < n$. Moreover, let $D$ be a Dyck path of type $B_n$. As in Section~\ref{sec:DyckCox} for usual Dyck paths, we identify $D$ with the set $\{b_{ij}\}$ of boxes below $D$. Label every box $b_{ij}$ with $j < n$ by $s_{n-1-i}$ and $b_{ij}$ with $j \geq n$ by $s_{2(n-1)-(i+j)}$. The bijection between Dyck paths in $\D_{B_n}$ and sortable elements is defined by mapping $D \in \D_{B_n}$ to $\psi_{B_n}(D) := \prod s$ where the product ranges all simple transpositions in the boxes $b_{ij}$ in the order as indicated in Figure~\ref{fig:dyckB5Cox}. For example, the Dyck path shown in the figure is mapped to
      $$s_5s_4s_3s_2s_1s_0 | s_5s_4s_2s_1s_0 | s_5s_2s_1 = (2,-2)(3,-6,-3,6)(4,5) = [1,-2,-6,5,4,3] \in \Cox_{B_6}.$$
      \begin{figure}
        \centering
          \begin{tikzpicture}[scale=1]
            \boxcollection{0.5}{(0,0)}{
            -2/1/,-2/2/,-2/3/,
            -1/1/,-1/2/,-1/3/,-1/4/,
            0/1/,0/2/,0/3/,0/4/,0/5/,
            1/1/,1/2/,1/3/,1/4/,1/5/,
            2/1/,2/2/,2/3/,2/4/,
            3/1/,3/2/,3/3/,
            4/1/,4/2/,
            5/1/}
            \latticepath{0.5}{(0,-3)}{very thick}{
              0/4/black,3/0/black,0/4/black,1/0/black}
            \latticepath{0.5}{(0,-3)}{very thick}{
              0/4/black,3/0/black,0/2/black,3/0/black}
            \latticepath{0.5}{(0,-3)}{very thin}{
              6/6/black,-3/3/black}
            \latticepath{0.5}{(0,-2.5)}{very thin}{
              5.5/5.5/grey}
            \latticepath{0.5}{(0,-2)}{very thin}{
              4.5/4.5/grey,0/1/grey}
            \latticepath{0.5}{(0,-1.5)}{very thin}{
              3.5/3.5/grey,0/2/grey}
            \latticepath{0.5}{(0,-1)}{very thin}{
              2.5/2.5/grey,0/2.5/grey}
            \latticepath{0.5}{(0,-0.5)}{very thin}{
              1.5/1.5/grey,0/2.5/grey}
            \latticepath{0.5}{(0,0)}{very thin}{
              0.5/0.5/grey,0/2.5/grey}
            \content{0.5}{(0.5,-2.5)}{
              0/0/$\color{black}s_5$,-1/1/$\color{black}s_4$,-2/2/$\color{black}s_3$,-3/3/$\color{black}s_2$,-4/4/$\color{black}s_1$,-5/5/$\color{black}s_0$,
              -1/0/$\color{black}s_5$,-2/1/$\color{black}s_4$,-4/3/$\color{black}s_2$,-5/4/$\color{black}s_1$,
              -2/0/$\color{black}s_5$,-5/3/$\color{black}s_2$,-6/4/$\color{black}s_0$,
              -6/3/$\color{black}s_1$}
            \content{0.5}{(-0.5,-1.5)}{
              -1/0/$\scriptstyle 6$,0/0/$\scriptstyle 5$,1/0/$\scriptstyle 4$,2/0/$\scriptstyle 3$,3/0/$\scriptstyle 2$,4/0/$\scriptstyle 1$}
            \latticepath{0.5}{(0,1.5)}{very thick}{
              4/0/white}
            \latticepath{0.5}{(0,1.47)}{very thick}{
              4/0/white}
            \latticepath{0.5}{(0,1.3)}{very thick}{
              4/0/white}
            \latticepath{0.5}{(0,1.27)}{very thick}{
              4/0/white}
            \latticepath{0.5}{(0,1.1)}{very thick}{
              4/0/white}
            \latticepath{0.5}{(0,1.13)}{very thick}{
              4/0/white}
            \draw [->] (-0.5,-3.0) -- (-0.2,-2.7);
            \draw [->] (-0.5,-2.5) -- (-0.2,-2.2);
            \draw [->] (-0.5,-2.0) -- (-0.2,-1.7);
            \draw [->] (-0.5,-1.5) -- (-0.2,-1.2);
            \draw [->] (-0.5,-1.0) -- (-0.2,-0.7);
            \draw [->] (-0.5,-0.5) -- (-0.2,-0.2);
          \end{tikzpicture}
        \caption{A Dyck path $\D_{B_6}$ with boxes labeled by simple transpositions.}
        \label{fig:dyckB5Cox}
      \end{figure}
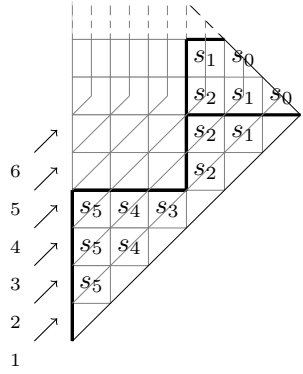
    To see that $\sigma = \psi_{B_n}(D) = X_1 X_2 \cdots X_n$ is in fact in $\Cox_{B_n}$, we have to show that it is a reduced expression for~$\sigma$, as the inclusion property $X_1 \supseteq X_2 \supseteq \ldots \supseteq X_n$ is given by construction.
    \begin{proposition}
      $X_1 X_2 \cdots X_n$ is a reduced expression for~$\sigma$.
    \end{proposition}
    \begin{proof}
      If $s_i s_{i-1}$ occurs in~$X_j$ and in~$X_{j+1}$ for some~$i$ and~$j$ then $s_{i-2}$ occurs also in~$X_j$ except for the case $i = 1$. But if $s_1s_0$ occurs in~$X_j$ and in~$X_{j+1}$ and~$s_1$ also occurs in~$X_{j+2}$ then~$s_2$ occurs in~$X_{j+2}$ left of~$s_1$. Thus, $X_1 X_2 \cdots X_n$ is reduced.
    \end{proof}
    This proposition immediately implies Theorem~\ref{th:bijectionsB}\eqref{eq:areatoinvB} for $\psi_{B_n}:\D_{B_n}\tilde\longrightarrow\NC_{B_n}$.
    \begin{corollary}
      Let $D \in \D_{B_n}$ and let $\sigma = \psi_{B_n}(D)$. Then
      $$\area(D) = \inv_B(\sigma).$$
    \end{corollary}

    To prove Theorem~\ref{th:bijectionsB}\eqref{eq:majtomajB} for $\psi_{B_n}:\D_{B_n}\tilde\longrightarrow\NC_{B_n}$, we use the fact that a Dyck path $D$ in $\D_{B_n}$ consists of a \lq\lq lower part\rq\rq\ $D'$ which is a Dyck path in $\Dn$, and an \lq\lq upper part\rq\rq\ $D''$. The path $D'$ is obtained from $D$ by replacing all north steps after the~$n$-th north step by east steps and $D''$ is obtained as the suffix of $D$ after the~$n$-th north step. For example, the Dyck path in $\D_{B_6}$ in Figure~\ref{fig:dyckB5Cox} consists of a lower part ending in $3$ east steps which is the Dyck path shown in Figure~\ref{fig:dyck5Cox}, and an upper part given by the suffix $NNE$.

    As $s_i$ and $s_j$ commute for $|i-j| > 1$, we can write $\psi_{B_n}(D)$ as $\psi_n(D')$ multiplied with the product of the boxes below $D''$ row by row from right to left and from bottom to top. Set $\sigma, \sigma'$ and $\sigma''$ to be the signed permutations associated to $D,D'$ and $D''$. For example,
    \begin{align*}
      \sigma = \psi_{B_n}(D)  &= s_5s_4s_3s_2s_1s_0 | s_5s_4s_2s_1s_0 | s_5s_2s_1 \\
                              &= s_5s_4s_3s_2s_1 | s_5s_4s_2 | s_5 \cdot s_0s_1s_2|s_0s_1 \\
                              &= \sigma' \hspace{5pt} \cdot \hspace{5pt} \sigma''.
    \end{align*}
    To have the example at hand for the following steps, we have
    $$
      \fmaj(D) = 2(4+5) = 18, \quad \fmaj(D') = 2 (6+5) = 22, \quad \fmaj(D'') = 0,
    $$
    where, for simplicity, we set $\fmaj(D'') := \sum_{i \in \Des(D'')}2(k-i)$ with $k$ being the number of steps in $D''$. Moreover, get have $\sigma = [1,-2,-6,5,4,3], \sigma' = [6,2,1,5,4,3], \sigma'' = [3,-2,-1,4,5,6]$, and thus,
    \begin{align*}
      \fmaj(\sigma) = 2(1+2+4+5)+2 = 26, &\quad \fmaj(\sigma') = 2(1+2+4+5) = 24, \quad \fmaj(\sigma'') = 2 \cdot 1 = 2 \\
      \ifmaj(\sigma) = 2(1+3+4+5)+2 = 28, &\hspace*{7.3pt} \ifmaj(\sigma') = 2(1+3+4+5) = 26, \hspace*{7.3pt} \ifmaj(\sigma'') = 2 \cdot 1 = 2.
    \end{align*}
    As we have seen in Section~\ref{sec:DyckCox}, we have, when considering $D'$ in $\Dn$, that $\maj(D') + \maj(\sigma') + \imaj(\sigma') = n(n-1)$. This gives, when considered in $\D_{B_n}$,
    \begin{align*}
      \fmaj(D') + \fmaj(\sigma') + \ifmaj(\sigma') = 2 n(n-1) + 2n = 2n^2.
    \end{align*}
    Using this fact, we show in three steps that
    \begin{align*}
      \fmaj(D) + \fmaj(\sigma) + \ifmaj(\sigma) &= \fmaj(D') + \fmaj(\sigma') + \ifmaj(\sigma').
    \end{align*}
    Theorem~\ref{th:bijectionsB}\eqref{eq:majtomajB} for $\psi_{B_n}:\D_{B_n}\tilde\longrightarrow\Cox_{B_n}$ then follows.
    \begin{lemma}\label{le:fmajD}
      Let $D,D',D'',$ and~$\sigma$ as above. Then $\fmaj(D) = \fmaj(D') + \fmaj(D'') - 2\neg(\sigma)$.
    \end{lemma}
    \begin{proof}
      By definition, $\fmaj(D) = \fmaj(D')+\fmaj(D'')-2(n-\neg(D))$. Observe here that $D''$ always starts with an east step, and thus we do not cut $D$ here at a descent. As $\neg(\sigma)$ is given by the number of $s_0$ occurring in the word for~$\sigma$ as given above, we obtain $\neg(\sigma) = n-\neg(D)$. The lemma follows.
    \end{proof}
    \begin{lemma}\label{le:fmajsigma}
      Let $\sigma,\sigma',$ and $D''$ as above. Then $\fmaj(\sigma) = \fmaj(\sigma') - \fmaj(D'') + \neg(\sigma)$.
    \end{lemma}
    \begin{proof}
      We are going to show that $\Des(\sigma') = \Des(\sigma) \uplus \{ k-i : i \in \Des(D'') \}$. The lemma can then be deduced using the definition of the flag major index. First, we observe that $\sigma_\ell = \sigma'_\ell$ for $\ell > k$, which gives that $\ell>k$ is a descent of~$\sigma$ if and only if it is a descent of~$\sigma'$. From Proposition~\ref{pr:lengthoflastdescent}, we obtain that $\{1,\ldots,k-1\} \subseteq \Des(\sigma')$ and that $k \notin \Des(\sigma')$. As $\sigma_k \leq 1$, $k$ is not a descent of~$\sigma$ either. Thus, it is left to show that
      $$ \{ i<k : i \in \Des(\sigma) \} \uplus \{ k-i : i \in \Des(D'') \} = \{1,\ldots,k-1\}.$$
      To see this, we will explicitly describe both descent sets using $\sigma''$. By construction,
      $$\Neg(\sigma'') = \big\{i+1 : s_i \text{ is the rightmost simple transposition in } X_\ell \text{ for some } \ell \big\},$$
      and the images of $\Neg(\sigma'')$ under $\sigma''$ are the negatives of the first $\neg(\sigma'')$ integers in increasing order and the image of the complement of $\Neg(\sigma'')$ are the last $k-\neg(\sigma'')$ integers also in increasing order. E.g., for $\sigma'' = s_0s_1s_2 | s_0s_1$ as above, we get $\Neg(\sigma'') = \{3,2\}$, and
      $$\sigma''(3) = -1, \quad \sigma''(2) = -2, \quad \sigma''(1) = 3, \sigma''(4) = 4.$$
      Using this description, we finally get
      \begin{align*}
        \{ i<k : i \in \Des(\sigma) \}  &= \{ i<k : i \notin \Neg(\sigma'') \hspace*{3.5pt}\text{ or  } \hspace*{3.5pt} i+1 \in \Neg(\sigma'') \},\\
        \{ k-i : i \in \Des(D'') \}      &= \{ i<k : i \in \Neg(\sigma'') \text{ and } i+1 \notin \Neg(\sigma'') \}.
      \end{align*}
      This completes the proof.
    \end{proof}
    \begin{lemma}\label{le:ifmajsigma}
      Let~$\sigma$ and~$\sigma'$ as above. Then $\ifmaj(\sigma) = \ifmaj(\sigma') + \neg(\sigma)$.
    \end{lemma}
    \begin{proof}
      The way $\sigma''$ is constructed, we obtain that $(\sigma'')^{-1}_i < (\sigma'')^{-1}_{i+1}$ and thus, $\iDes(\sigma'') = \emptyset$. By definition, $\sigma^{-1} = (\sigma'')^{-1} (\sigma')^{-1}$, which gives $\iDes(\sigma) = \iDes(\sigma')$. The lemma follows with the fact that $\neg(\sigma') = 0$.
    \end{proof}

  \bibliographystyle{amsplain} 
  \bibliography{../../mrabbrev,../../bibliography}
\end{document}